# Hellinger vs. Kullback-Leibler multivariable spectrum approximation

Augusto Ferrante, Michele Pavon and Federico Ramponi

**Abstract**

In this paper, we study a matricial version of the Byrnes-Georgiou-Lindquist generalized moment problem with complexity constraint. We introduce a new metric on multivariable spectral densities induced by the family of their spectral factors which, in the scalar case, reduces to the Hellinger distance. We solve the corresponding constrained optimization problem via duality theory. A highly nontrivial existence theorem for the dual problem is established in the Byrnes-Lindquist spirit. A matricial Newton-type algorithm is finally provided for the numerical solution of the dual problem. Simulation indicates that the algorithm performs effectively and reliably.

## I. INTRODUCTION

In the past ten years, building on their previous work, Byrnes, Georgiou, Lindquist and collaborators have developed a broad program on generalized analytic interpolation and generalized moment problems that arise in spectral estimation and robust control [9], [6], [24], [7], [10], [18], [22], [23], [25], [26], [42], [3], [27], [28], [12], [8]. While we refer the reader to the cited literature for better motivation, we recall that many problems of $H^\infty$ control, signal processing and maximal power transfer in circuit theory may be reduced to a Nevanlinna-Pick interpolation

Work partially supported by the MIUR-PRIN Italian grant "Identification and Control of Industrial Systems".

A. Ferrante and F. Ramponi are with the Dipartimento di Ingegneria dell'Informazione, Università di Padova, via Gradenigo 6/B, 35131 Padova, Italy augusto@dei.unipd.it, rampo@dei.unipd.it

M. Pavon is with the Dipartimento di Matematica Pura e Applicata, Università di Padova, via Trieste 63, 35131 Padova, Italy pavon@math.unipd.it

problem, see e.g. [17], [6], [51]. In all of these applications, it is crucial to put a bound on the degree of the interpolant so that the controller, filter, etc. has limited complexity. As is well known, while the Nevanlinna-Pick theory features a simple criterion in terms of the Pick matrix for the existence of solutions and beautiful iterative techniques (Schur-type algorithms) to produce solutions when they exist, the degree specification on the interpolant is much harder to capture in this framework. The overcoming of this difficulty by the Byrnes, Georgiou, Lindquist school has open the way to several new applications in speech processing, bioengineering and robust control, see [5], [43], [33]. Notice that [26], [3], [28] deal with the more difficult multidimensional case.

One of the central steps, in these authors approach, is the formulation of a convex optimization problem that includes as a (very) special case maximum entropy problems. The smooth parametrization of the complete class of interpolants occurs in the optimization setting, where it is crucial the dependence of the criterion on certain "a priori" parameters, cf. e.g. Remark 3.2 below. It should be observed that, as in all of the previously mentioned applications, the primal problem is infinite dimensional, while the dual problem is finite dimensional. Hence, it is natural to seek the (unique) solution to the primal problem via duality theory.

In [32], Georgiou and Lindquist, have applied this convex optimization approach to constrained spectrum approximation. The basic ingredients of the optimization problem are the following: An a priori power spectral density $\Psi$ is given. Then new data become available in the form of asymptotic state-covariance statistics for a bank of filters. The latter induces a linear constraint on the family of spectral densities. It is then necessary to find a spectrum $\Phi$ that satisfies the constraint and is as close as possible to $\Psi$ in a prescribed metric.

In [32], a Kullback-Leibler criterion was employed, where minimization is performed with respect to the second argument. This unusual choice was dictated by two considerations: i) The desire to have maximum entropy as a special case ($\Psi = I$); ii) the simple form of the optimal solution belonging to a parametric family of "rational" densities. The latter class, as well as another parametric class of "exponential type" [27, p.3], were recognized from the start [9], [7] to be critical points of logarithmic entropy-like functionals. In [27] and in [28], homotopy like methods were proposed as an effective tool to solve a class of scalar and multidimensional generalized moment problems.

In this paper, we seek to investigate constrained approximation of spectral density functions in a different metric, also originating in mathematical statistics, namely the *Hellinger distance* [15], [39], [36], [40]. Observing that this distance amounts to the minimum $L^2$ distance between spectral factors of the two spectra, we are led to a natural extension of this metric to multivariable spectra (Theorem 6.1). We then show that the latter observation opens the way to attack the general multivariable case where, differently from the Kullback-Leibler case, an explicit form for the optimal solution may be obtained, see Theorem 7.2 . We also establish an existence theorem for the dual problem (Theorem 7.6) that paralles a correspondig fundamental result due to Byrnes and Lindquist [12]. We finally investigate iterative numerical methods to solve the dual problem. Although the dual problem is an unconstrained convex finite dimensional problem, the numerics is nontrivial. As observed in [3, Section VI], [32, SectionVII], the dual functional has unbounded gradient at the boundary. Reformulation of the problem to avoid this difficulty may lead to loss of global convexity, requiring to initialize any descent method close to the minimum [6], [18], [42], [3]. As in [45] for the Kullback-Leibler case, we prefer to employ matricial descent methods. A number of nontrivial difficulties in the algorithms are overcome by resorting to ideas and results from spectral factorization theory. In our simulation, these iterative schemes (particularly a Newton type method with backtracking line search) appear to perform effectively and reliably.

The paper is outlined as follows. Section II is devoted to the formulation of a generalized moment problem in the sense of Byrnes-Georgiou-Lindquist, and to the corresponding existence question. In Section III, two approximation problems for scalar spectral densities are introduced. The first employs a Kullback-Leibler type criterion, the second features the Hellinger distance. Optimality conditions for these two problems are presented in Section IV. The multivariable version of the two approximation problems, and the corresponding difficulties in the variational analysis, are discussed in Section V. Section VI is devoted to the introduction of a new metric on multivariable spectral densities induced by the corresponding spectral factors. The multivariable spectrum approximation problem with respect to the distance of Section VI is solved in Section VII. Finally, Section VIII deals with the numerical solution of the dual problem.

## II. A GENERALIZED MOMENT PROBLEM

We consider the following basic set-up patterned after [26], [32], [28]. Let $\mathcal{S}_+^{m \times m}(\mathbb{T})$ be the family of bounded, coercive, $\mathbb{C}^{m \times m}$-valued spectral density functions on the unit circle. Thus, a measurable, bounded matrix valued function $\Phi$ belongs to $\mathcal{S}_+^{m \times m}(\mathbb{T})$ if it satisfies the following properties:

- the values of $\Phi$ are $m \times m$, Hermitian, non negative definite matrices;
- there exists a positive constant $c_\Phi$ such that $\Phi(e^{i\theta}) - c_\Phi I$ is positive definite a.e. on $\mathbb{T}$.

Notice that $\Phi \in \mathcal{S}_+^{m \times m}(\mathbb{T})$ if and only if $\Phi^{-1} \in \mathcal{S}_+^{m \times m}(\mathbb{T})$. Let $\Psi \in \mathcal{S}_+^{m \times m}(\mathbb{T})$ represent an *a priori* estimate of the spectrum of an underlying zero-mean, wide-sense stationary $m$-dimensional stochastic process $\{y(n), n \in \mathbb{Z}\}$. We consider a rational transfer function

$$G(z) = (zI - A)^{-1}B, \qquad A \in \mathbb{C}^{n \times n}, B \in \mathbb{C}^{n \times m}, \tag{1}$$

where $A$ has all its eigenvalues in the open unit disc, $B$ is full column rank, and $(A, B)$ is a reachable pair. Here $G$ models a *bank of filters*. We consider the situation where new data become available in the form of an asymptotic estimate $\Sigma > 0$ of the state covariance of the system with transfer function $G$ and input the unknown process $y$. In other words, we suppose we can estimate the covariance of the $n$- dimensional stationary process $\{x_k; k \in \mathbb{Z}\}$ satisfying

$$x_{k+1} = Ax_k + By_k, \quad k \in \mathbb{Z}. \tag{2}$$

In general, $\Psi$ is not consistent with $\Sigma$, and it is necessary to find $\Phi$ in $\mathcal{S}_+^{m \times m}(\mathbb{T})$ that is closest to $\Psi$ in a suitable sense among spectra consistent with $\Sigma$, namely satisfying

$$\int G\Phi G^* = \Sigma, \tag{3}$$

where star denotes transposition plus conjugation. Here, and throughout the paper, integration takes place on $[-\pi, \pi]$ with respect to the normalized Lebesgue measure $d\theta/2\pi$. The question of existence of $\Phi \in \mathcal{S}_+^{m \times m}(\mathbb{T})$ satisfying (3) and, when existence is granted, the parametrization of all solutions to (3), may be viewed as a generalized moment problem. For instance, in the case $m = 1$, take $G(z)$ with $k$-th component $G_k(z) = z^{k-n-1}$. Take moreover

$$A = \begin{bmatrix} 0 & 1 & 0 & \ldots & 0 \\ 0 & 0 & 1 & \ldots & 0 \\ \vdots & \vdots & & \ddots & \vdots \\ 0 & 0 & 0 & \ldots & 1 \\ 0 & 0 & 0 & \ldots & 0 \end{bmatrix}, \quad B = \begin{bmatrix} 0 \\ 0 \\ \vdots \\ 0 \\ 1 \end{bmatrix} \tag{4}$$

Here $c_k := E\{y(n)\bar{y}(n+k)\}$. This is the *covariance extension problem*, where the information available on the process $y$ is the finite sequence of covariance lags $c_0, c_1, \ldots, c_{n-1}$. It is known that the set of densities consistent with the data is nonempty if $\Sigma \geq 0$ and contains infinitely many elements if $\Sigma > 0$ [34], see also [22], [9], [10], [23].

*Existence* of $\Phi \in \mathcal{S}_+^{m \times m}(\mathbb{T})$ satisfying constraint (3) is a nontrivial issue. It was shown in [24], [25] that such family is nonempty if and only if there exists $H \in \mathbb{C}^{m \times n}$ such that

$$\Sigma - A\Sigma A^* = BH + H^*B^*,\tag{5}$$

or, equivalently, the following rank condition holds

$$\operatorname{rank} \begin{pmatrix} \Sigma - A\Sigma A^* & B \\ B^* & 0 \end{pmatrix} = 2m.\tag{6}$$

We wish to give an alternative formulation of this existence result. First of all, notice that W.L.O.G. we can take $\Sigma = I$. Indeed, if $\Sigma \neq I$, it suffices to replace $G$ with $G' := \Sigma^{-1/2}G$ and $(A, B)$ with $(A' = \Sigma^{-1/2}A\Sigma^{1/2}, B' = \Sigma^{-1/2}B)$. Thus constraint (3) from now on reads

$$\int G\Phi G^* = I.\tag{7}$$

Let $\Pi_B = B(B^*B)^{-1}B^*$ denote the orthogonal projection onto $\operatorname{Range}(B)$.

*Proposition 2.1:* A necessary and sufficient condition for the existence of spectra in $\mathcal{S}_+^{m \times m}(\mathbb{T})$ satisfying (7) is that the following relation holds

$$(I - \Pi_B)(I - AA^*)(I - \Pi_B) = 0.\tag{8}$$

When (8) is satisfied, there exists $\Phi \in \mathcal{S}_+^{m \times m}(\mathbb{T})$ satisfying (7) of McMillan degree less than or equal to $2n$.

*Proof: Necessity:* Suppose there exists $y$ $m$-dimensional, wide-sense stationary with spectral density $\Phi \in \mathcal{S}_+^{m \times m}(\mathbb{T})$ satisfying (7). Let $x$ be defined by (2). Taking covariances on both sides of (2), we get

$$I = AA^* + AE\{x_k y_k^*\}B^* + BE\{y_k x_k^*\}A^* + BE\{y_k y_k^*\}B^*.$$

The latter relation implies (8).

*Sufficiency:* We adapt the argument in [26, p.1814]. For a given purely non deterministic $m$-dimensional process $y$ with spectrum $\Phi$, define the process $w$ as the output of the linear stable system

$$x_{k+1} = Ax_k + By_k,\tag{9}$$

$$w_k = (B^*B)^{-1}B^*x_{k+1}\tag{10}$$

Inverting the system (9)-(10), we get

$$x_{k+1} = (I - \Pi_B) A x_k + B w_k, \tag{11}$$

$$y_k = -(B^*B)^{-1} B^* A x_{k-1} + w_k. \tag{12}$$

Write (8) as a Lyapunov identity

$$I = (I - \Pi_B) A A^* (I - \Pi_B) + \Pi_B. \tag{13}$$

Since $(A, B)$ is controllable, so is the pair $\left((I - \Pi_B)A, B(B^*B)^{-1/2}\right)$. It now follows from (13) that $(I - \Pi_B)A$ has all eigenvalues in the open unit disc $\mathbb{D}$. Thus, system (11)-(12) is stable, $(B^*B)^{-1}B^*G(z)$ is minimum phase and the processes $y$ and $w$ are causally equivalent. It follows that if we choose $w$ to be a white noise sequence with intensity $E\{w_k w_k^*\} = (B^*B)^{-1}$ and $y$ to be defined by (11)-(12) then: i) (11)-(12) is the innovation representation of $y$; ii) the state covariance of the steady-state Kalman filter (11)-(12) satisfies the Lyapunov equation (13) and is therefore the identity; iii) the spectral density of $y$ is given by

$$\Phi_y = W(z)(B^*B)^{-1}W(z)^*, \tag{14}$$

where

$$W(z) = I - (B^*B)^{-1} B^* A \left( zI - (I - \Pi_B)A \right)^{-1} B,$$

is the transfer function of (11)-(12). We conclude that if we feed $G$ in (9) with such a process $y$, the filter state $x$ will have the required covariance, namely the identity matrix, and (7) will be satisfied. Moreover, $\Phi_y$ is rational of McMillan degree at most $2n$ and it belongs to $\mathcal{S}_+^{m \times m}(\mathbb{T})$ since its values and the values of $\Phi_y^{-1}$ on $\mathbb{T}$ are positive definite matrices. The spectrum (14) has been shown in [26, Section III] to be the maximum entropy spectrum among those satisfying (7). This is there accomplished in a clever way, by relating the constrained maximum entropy problem to a special one-step-ahead prediction problem.

## III. Constrained spectrum approximation: The scalar case

### A. Kullback-Leibler criterion

In [32], the Kullback-Leibler measure of distance for spectra in $\mathcal{S}_+(\mathbb{T}) := \mathcal{S}_+^{1 \times 1}(\mathbb{T})$ was introduced:

$$\mathbb{D}(\Psi \| \Phi) = \int \Psi \log \left( \frac{\Psi}{\Phi} \right).$$

As is well known, this pseudo-distance originates in hypothesis testing, where it represents the mean information for observation for discrimination of an underlying probability density from another [38, p.6]. It also plays a central role in information theory, identification, stochastic processes, etc., see e.g. [37], [16], [14], [21], [2], [13], [49], [46] and references therein. It is also known in these fields as *divergence*, *relative entropy*, *information distance*. etc. If

$$\int \Phi = \int \Psi,$$

we have $\mathbb{D}(\Psi \| \Phi) \geq 0$. The choice of $\mathbb{D}(\Psi \| \Phi)$ as a distance measure, even for spectra that have different zeroth moment, is discussed in [32, Section III]. It is observed there that the constraint (3) often fixes the zeroth Fourier coefficient of feasible spectra (this happens for sure when $A$ is singular). In that case, rescaling $\Psi$, we are guaranteed that the index is nonnegative and equal to zero if and only if the two spectra are equal. T. Georgiou has kindly informed us [29], that even when $A$ is nonsingular, under a rather mild assumption, it is possible to modify the index so that all $\Phi$ satisfying the constraint have the same zeroth moment. In any case, the method entails a rescaling of the a priori density $\Psi$, so that the optimization problem amounts to approximating the "shape" of the a priori spectrum. This is of course sensible to pursue in several engineering applications such as speech processing.

We mention, for the benefit of the reader, that, in the same spirit, Georgiou has very recently investigated other distances for power spectra, [30], [31]. Motivated by classical prediction theory, where the optimal one step ahead predictor does not depend on the $L^1$ norm of the spectrum, he seeks natural distances between *rays* of spectral densities. Considering the degradation of performance when an optimal predictor for one stochastic process is employed to predict a different stochastic process, he is naturally led to introduce a certain metric on rays.

As observed in the introduction, notice that minimizing $\mathbb{D}(\Psi \| \Phi)$ rather than $\mathbb{D}(\Phi \| \Psi)$ is unusual with respect to the statistics-probability-information theory world. Besides leading to a more tractable form of the optimal solution, however, it also includes as special case ($\Psi \equiv 1$) maximization of entropy [26]. In [32], the following problem is considered:

*Problem 3.1:* (**Approximation problem 1**) Given $\Psi \in \mathcal{S}(\mathbb{T})$, find $\hat{\Phi}_{KL}$ that solves

$$\text{minimize} \quad \mathbb{D}(\Psi \| \Phi) \tag{15}$$

$$\text{over} \quad \left\{ \Phi \in \mathcal{S}(\mathbb{T}) \mid \int G \Phi G^* = I \right\}. \tag{16}$$

*Remark 3.2:* In the context of the covariance extension problem (4), the minimizers of (3.1), when $\Psi$ ranges over positive trigonometric polynomials of degree $n$, are precisely the coercive spectra consistent with the first $n$ covariance lags and of degree at most $2n$, [9], [10], [22], [23]. This illustrates the role of the "a priori parameter" $\Psi$ in obtaining a description of all solution to the moment problem of prescribed complexity.

## B. Hellinger criterion

In this paper, we consider a different metric on spectral density functions. Given $\Phi$ and $\Psi$ in $\mathcal{S}(\mathbb{T})$, the *Hellinger distance* is defined by

$$d_H(\Phi, \Psi) := \left[ \int_{-\pi}^{\pi} \left( \sqrt{\Phi(e^{i\theta})} - \sqrt{\Psi(e^{i\theta})} \right)^2 \frac{d\theta}{2\pi} \right]^{1/2}.$$

It is a *bona fide* distance on $\mathcal{S}(\mathbb{T})$. Moreover, it satisfies the following properties.

*Proposition 3.3: Consider* $\Phi, \Psi \in \mathcal{S}(\mathbb{T})$. *Then*

1) $d_H(\Phi, \Psi) \leq \sqrt{\|\Phi\|_1 + \|\Psi\|_1}$;
2) $d_H(\Phi, \Psi)^2 \leq \|\Phi - \Psi\|_1$;
3) $\|\Phi - \Psi\|_1 \leq \left( \sqrt{\|\Phi\|_1} + \sqrt{\|\Psi\|_1} \right) d_H(\Phi, \Psi)$.

These extend well-known properties of the Hellinger distance in the case of probability density functions. The straightforward proof may be found in [20].

*Remark.* On a finite-dimensional statistical manifold, endowed with the Fisher information as the metric tensor, both the Hellinger distance and the Kullback-Leibler pseudo-distance can be viewed as instances of the broader concept of $\alpha$-*divergences* between two points, which arise from the so-called Amari connections. In particular, the $0$-divergence, which indeed is the Hellinger distance, arises from the Levi-Civita connection. See [1, p. 66 and following].

We consider the following approximation problem.

*Problem 3.4:* (**Approximation problem 2**) Given $\Psi \in \mathcal{S}(\mathbb{T})$, find $\hat{\Phi}_H$ that solves

$$\text{minimize} \quad d_H^2(\Phi, \Psi) \tag{17}$$

$$\text{over} \quad \left\{ \Phi \in \mathcal{S}(\mathbb{T}) \mid \int G \Phi G^* = I \right\}. \tag{18}$$

# IV. OPTIMALITY CONDITIONS

## A. Kullback-Leibler approximation

Consider first Problem 3.1. The variational analysis in [32] is outlined as follows (see also [45]). For $\Lambda \in \mathbb{C}^{n \times n}$ Hermitian satisfying $G^* \Lambda G > 0$ on all of $\mathbb{T}$, consider the *Lagrangian function*

$$
\begin{aligned}
L(\Phi, \Lambda) &= \mathbb{D}(\Psi \| \Phi) + \text{tr}\left( \Lambda \left( \int G \Phi G^* - I \right) \right) \\
&= \mathbb{D}(\Psi \| \Phi) + \int G^* \Lambda G \Phi - \text{tr}(\Lambda), \quad (19)
\end{aligned}
$$

where "$\text{tr}$" denotes the trace operator. Consider the *unconstrained* minimization of the strictly convex functional $L(\Phi, \Lambda)$:

$$
\text{minimize}\{L(\Phi, \Lambda) | \Phi \in \mathcal{S}(\mathbb{T})\} \quad (20)
$$

This is a convex optimization problem. The variational analysis yields the following result.

*Theorem 4.1:* The unique solution $\hat{\Phi}_{KL}$ to problem (20) is given by

$$
\hat{\Phi}_{KL} = \frac{\Psi}{G^* \Lambda G}. \quad (21)
$$

*Moreover, suppose $\hat{\Lambda} = \hat{\Lambda}^*$ is such that*

$$
G^* \hat{\Lambda} G > 0, \quad \forall e^{i\theta} \in \mathbb{T}, \quad (22)
$$

$$
\int G \frac{\Psi}{G^* \hat{\Lambda} G} G^* = I. \quad (23)
$$

*Then $\hat{\Phi}_{KL}$ given by*

$$
\hat{\Phi}_{KL} = \frac{\Psi}{G^* \hat{\Lambda} G} \quad (24)
$$

*is the unique solution of the approximation Problem (3.1).*

Thus, the original Problem 3.1 is now reduced to finding $\hat{\Lambda}$ satisfying (22)-(23). This is accomplished via duality theory. Consider the dual functional

$$
\Lambda \to \inf\{L(\Phi, \Lambda) | \Phi \in \mathcal{S}(\mathbb{T})\}.
$$

For $\Lambda$ satisfying (22), the dual functional takes the form

$$
\Lambda \to L(\frac{\Psi}{G^* \Lambda G}, \Lambda) = \int \Psi \log G^* \Lambda G - \text{tr}(\Lambda) + \int \Psi. \quad (25)
$$

Consider now the maximization of the dual functional (25) over the set

$$
\mathcal{L}^{KL} := \{\Lambda = \Lambda^* | G^* \Lambda G > 0, \forall e^{i\theta} \in \mathbb{T}\}. \quad (26)
$$

Let, as in [32],

$$J_\Psi(\Lambda) := - \int \Psi \log G^* \Lambda G + \mathrm{tr}\,(\Lambda)$$

The dual problem is then equivalent to

$$\text{minimize} \quad \{J_\Psi(\Lambda)|\Lambda \in \mathcal{L}^{KL}\}. \tag{27}$$

The dual problem is also a convex optimization problem. In [32], $\Lambda$ is further restricted to belong to the range of the operator $\Gamma$ defined on the set $\mathcal{C}_H(\mathbb{T})$ of Hermitian-valued continuos functions defined in $\mathbb{T}$, by

$$\Gamma(\Phi) = \int G \Phi G^*, \qquad \Phi \in \mathcal{C}_H(\mathbb{T}). \tag{28}$$

As mentioned in Section II (equation (5)),

$$\mathrm{Range}(\Gamma) = \{\Sigma = \Sigma^* : \ \exists H \in \mathbb{C}^{m \times n} \text{ s.t. } \Sigma - A\Sigma A^* = BH + H^* B^*, \} \tag{29}$$

The problem then becomes

$$\text{minimize} \quad \{J_\Psi(\Lambda)|\Lambda \in \mathcal{L}_\Gamma^{KL}\}, \quad \mathcal{L}_\Gamma^{KL} = \mathcal{L}^{KL} \cap \mathrm{Range}(\Gamma). \tag{30}$$

The reason is that the orthogonal complement of $\mathrm{Range}(\Gamma)$ is given by

$$\mathrm{Range}(\Gamma)^\perp = \{M = M^* | G^* M G \equiv 0 \text{ on } \mathbb{T}\}. \tag{31}$$

This follows from the fact that $M \in \mathrm{Range}(\Gamma)^\perp$ iff $\forall \Phi \in \mathcal{C}_H(\mathbb{T})$

$$0 = \mathrm{tr}\,(\int G \Phi G^* M) = \int G^* M G \Phi.$$

The dual functional is shown in [32] to be strictly convex on the restricted domain $\mathcal{L}_\Gamma^{KL}$. It is also shown in [12] that $J_\Psi$ has a unique minimum point in $\mathcal{L}_\Gamma^{KL}$. This result implies that, under assumption (6), there exists a (unique) $\hat{\Lambda}$ in $\mathcal{L}_\Gamma^{KL}$ satisfying (23). Such a $\hat{\Lambda}$ then provides the optimal solution of the primal problem (15)-(16) via (24).

## B. Hellinger approximation

The variational analysis for Problem 3.4 is very similar, see [20] for details. We state without proof the following result: It will be proven in Section VII in the (more general) multivariable case.

*Theorem 4.2:* Assume that Problem 3.4 is feasible, namely that condition (6) (or, equivalently, condition (8)) is satisfied. Then there exists $\hat{\Lambda} = \hat{\Lambda}^* \in \mathbb{C}^{n \times n}$ such that

$$1 + G^* \hat{\Lambda} G > 0, \quad \forall e^{i\theta} \in \mathbb{T}, \qquad \int G \frac{\Psi}{(1 + G^* \hat{\Lambda} G)^2} G^* = I. \tag{32}$$

*In this case, Problem 3.4 admits a unique solution which is given by*

$$\hat{\Phi}_H = \frac{\Psi}{(1 + G^* \hat{\Lambda} G)^2} \tag{33}$$

*Remark 4.3:* Suppose the a priori density $\Psi$ is rational. Then, the solution in (33) has in general degree $2n$ higher than the solution in (24).

## V. Constrained spectrum approximation: The multivariate case

### A. Kullback-Leibler approximation

Multivariable Kullback-Leibler approximation has been investigated in [26], [28], whereas [3] deals with the multivariate Nevanlinna-Pick problem. In statistical quantum mechanics, the state of an $n$-level system is represented by a *density matrix* $\rho$, namely a Hermitian, positive-semidefinite matrix in $\mathbb{C}^{n \times n}$ with unit trace [48]. The convex set of density matrices has as extreme points the one dimensional projections. The latter can be identified with the *pure states* of the system $|\psi\rangle$, where $\psi$ is a unit vector in $\mathbb{C}^n$, via $\rho = \langle \psi, \cdot \rangle \psi$. Quantum analogues of entropy-like functionals have been considered since the early days of quantum mechanics [50]. Recently, renewed interest has originated in Quantum Information applications [44]. The quantum relative entropy between two density matrices is defined by:

$$\mathbb{D}(\rho||\sigma) := \mathrm{tr}\,(\rho(\log \rho - \log \sigma)). \tag{34}$$

Klein's inequality yields that $\mathbb{D}(\rho||\sigma) \geq 0$ if and only if $\rho = \sigma$. Moreover, as in the classical case, the quantum relative entropy is jointly convex in its arguments. We are then led to the following definition: Given $\Phi$ and $\Psi$ in $\mathcal{S}_+^{m \times m}(\mathbb{T})$, the relative entropy $\mathbb{D}(\Psi||\Phi)$ is given by

$$\mathbb{D}(\Psi||\Phi) = \int \mathrm{tr}\,(\Psi(\log \Psi - \log \Phi)). \tag{35}$$

First of all, we need to worry about nonnegativity of $\mathbb{D}(\Psi||\Phi)$ and whether it is zero iff $\Psi = \Phi$.

*Proposition 5.1:* Let $\Phi, \Psi \in \mathcal{S}_+^{m \times m}(\mathbb{T})$. Define $\Psi_1 = \Psi/\mathrm{tr}\,\Psi$ and $\Phi_1 = \Phi/\mathrm{tr}\,\Phi$. Then

$$\mathbb{D}(\Psi||\Phi) = \mathbb{D}(\mathrm{tr}\,\Psi||\mathrm{tr}\,\Phi) + \int (\mathrm{tr}\,\Psi)\mathrm{tr}\,(\Psi_1(\log \Psi_1 - \log \Phi_1)). \tag{36}$$

*It follows that when $\int \mathrm{tr}\,\Psi = \int \mathrm{tr}\,\Phi$, then $\mathbb{D}(\Psi\|\Phi) \geq 0$. Moreover, $\mathbb{D}(\Psi\|\Phi) = 0$ if and only if the two spectra coincide.*

*Proof:*

$$\mathbb{D}(\Psi\|\Phi) = \mathrm{tr}\int \Psi \left(\log\Psi - \log\Phi\right)$$

$$= \mathrm{tr}\int \mathrm{tr}\,(\Psi)\Psi_1 \left(\log \mathrm{tr}\,(\Psi)\Psi_1 - \log \mathrm{tr}\,(\Phi)\Phi_1\right)$$

$$= \mathrm{tr}\int \mathrm{tr}\,(\Psi)\Psi_1 \left((\log \mathrm{tr}\,(\Psi))I + \log\Psi_1 - (\log \mathrm{tr}\,(\Phi))I - \log\Phi_1\right)$$

$$= \mathrm{tr}\int \mathrm{tr}\,(\Psi)\Psi_1 \left(\log\Psi_1 - \log\Phi_1\right) + \mathrm{tr}\int \mathrm{tr}\,(\Psi)\Psi_1 \left((\log \mathrm{tr}\,(\Psi)) - (\log \mathrm{tr}\,(\Phi))\right)I$$

$$= \int \mathrm{tr}\,(\Psi)\mathrm{tr}\,\left(\Psi_1 \left(\log\Psi_1 - \log\Phi_1\right)\right) + \int \mathrm{tr}\,(\Psi_1)\mathrm{tr}\,\Psi \log\frac{\mathrm{tr}\,\Psi}{\mathrm{tr}\,\Phi}$$

$$= \int \mathrm{tr}\,(\Psi)\mathrm{tr}\,\left(\Psi_1 \left(\log\Psi_1 - \log\Phi_1\right)\right) + \mathbb{D}(\mathrm{tr}\,\Psi\|\mathrm{tr}\,\Phi).$$

Since $\mathrm{tr}\,\Psi_1(e^{i\theta}) = \mathrm{tr}\,\Phi_1(e^{i\theta}) = 1, \forall\theta \in [-\pi, \pi]$, it follows from Klein's inequality that

$$\mathrm{tr}\,\Psi_1(e^{i\theta}) \left(\log\Psi_1(e^{i\theta}) - \log\Phi_1(e^{i\theta})\right) \geq 0, \quad \forall\theta.$$

The latter implies that

$$\int (\mathrm{tr}\,\Psi)\mathrm{tr}\,\left(\Psi_1 \left(\log\Psi_1 - \log\Phi_1\right)\right) \geq 0.$$

When $\int \mathrm{tr}\,\Psi = \int \mathrm{tr}\,\Phi$, we also have $\mathbb{D}(\mathrm{tr}\,\Psi\|\mathrm{tr}\,\Phi) \geq 0$. Thus, when $\int \mathrm{tr}\,\Psi = \int \mathrm{tr}\,\Phi$, $\mathbb{D}(\Psi\|\Phi)$ is the sum of two nonnegative terms and the conclusion follows. $\square$

Consider again Problem 3.1.

*Problem 5.2:* (**Approximation problem 1**) For $\Psi \in \mathcal{S}_+^{m\times m}(\mathbb{T})$

$$\text{minimize} \quad \mathbb{D}(\Psi\|\Phi) \tag{37}$$

$$\text{over} \quad \left\{\Phi \in \mathcal{S}_+^{m\times m}(\mathbb{T}) \mid \int G\Phi G^* = I\right\}, \tag{38}$$

where $\mathbb{D}(\Psi\|\Phi)$ is defined by (35). As in the scalar case, an *a posteriori* rescaling of the prior density is in general necessary. In the light of Proposition 5.1, if $\hat\Phi$ is the solution of (5.2), the new prior is

$$\hat\Psi = \frac{\int \mathrm{tr}\,\hat\Phi}{\int \mathrm{tr}\,\Psi}\Psi.$$

For $\Lambda \in \mathbb{C}^{n \times n}$ Hermitian such that $G^* \Lambda G$ is positive definite on all of $\mathbb{T}$, define again the *Lagrangian*

$$
\begin{aligned}
L(\Phi, \Lambda) &= \mathbb{D}(\Psi \| \Phi) + \text{tr}\left(\Lambda\left(\int G\Phi G^* - I\right)\right) \\
&= \mathbb{D}(\Psi \| \Phi) + \int G^* \Lambda G \Phi - \text{tr}(\Lambda). \quad (39)
\end{aligned}
$$

The following step, entailing the *unconstrained* minimization of the strictly convex functional $L(\Phi, \Lambda)$ on $\Psi \in \mathcal{S}_+^{m \times m}(\mathbb{T})$, is a stumbling block. The optimality condition reads [28, Section IV]

$$
\int_0^\infty (\hat{\Phi}_{KL} + \tau I)\Psi(\hat{\Phi}_{KL} + \tau I)d\tau = G^* \Lambda G. \quad (40)
$$

In general, an explicit expression for $\hat{\Phi}_{KL}$ in terms of $\Psi$ and $\Lambda$ cannot be obtained, and the variational analysis ends here. We mention that the minimization with respect to the first argument of the relative entropy can instead be carried out explicitly, leading to a solution of the exponential form

$$
\Phi^o = c \exp(\log \Psi - G^* \Lambda G),
$$

see [28, Section IV]. Homotopy like methods are described in [28] to find $\Lambda$, when it exists, such that $\Phi^o$ satisfies the constraint.

## B. Hellinger approximation

Recall that, for a positive semidefinite Hermitian matrix $M$, $M^{1/2}$ is the *square root* of $M$, namely the unique Hermitian matrix whose square is $M$. If $V$ is a unitary matrix that diagonalizes $M$ so that $M = V^* \text{diag}(\alpha_1^2, \ldots, \alpha_m^2)V$, then simply $M^{1/2} = V^* \text{diag}(\alpha_1, \ldots, \alpha_m)V$. Motivated by the analogy with the Kullback-Leibler case, and by the scalar case, we define the Hellinger distance for $\Phi$ and $\Psi$ in $\mathcal{S}_+^{m \times m}(\mathbb{T})$ to be

$$
d_H^2(\Phi, \Psi) := \int_{-\pi}^{\pi} \text{tr}\left[\Phi^{1/2}(e^{i\theta}) - \Psi^{1/2}(e^{i\theta})\right]^2 \frac{d\theta}{2\pi}. \quad (41)
$$

Notice that (41) appears also as the natural generalization of the Hellinger distance for density operators of statistical quantum physics introduced in [41]. Consider again the strictly convex Problem 3.4:

$$
\text{minimize} \quad d_H^2(\Phi, \Psi) \quad (42)
$$

$$
\text{over} \quad \left\{\Phi \in \mathcal{S}_+^{m \times m}(\mathbb{T}) \mid \int G\Phi G^* = I\right\}, \quad (43)
$$

where $d_H^2(\Phi, \Psi)$ is now given by (41). Define $\mathcal{L}^H$ by

$$\mathcal{L}^H := \{\Lambda \in \mathbb{C}^{n \times n} | \Lambda = \Lambda^*, I + G^* \Lambda G > 0 \text{ a.e.on} \mathbb{T}\}. \tag{44}$$

For $\Lambda \in \mathcal{L}^H$, consider the *Lagrangian*

$$L^H(\Phi, \Lambda) = d_H^2(\Phi, \Psi) + \text{tr}\left(\Lambda\left(\int G\Phi G^* - I\right)\right).$$

The unconstrained minimization of the strictly convex functional $L^H$ over $\Phi \in \mathcal{S}_+^{m \times m}(\mathbb{T})$, however, leads to an optimality condition (expressing the the unique optimum $\hat{\Phi}_H$ in terms of $\Psi$ and $\Lambda$) which does not appear to be useful.

To obtain such optimality condition, we first need an expression for the directional derivative of the matrix square root. More precisely, given $P = P^* > 0$ let $S(P) := P^{1/2}$ and $\delta P = \delta P^*$: We want to compute

$$\delta S(P, \delta P) := \lim_{\varepsilon \to 0} \frac{(P + \varepsilon \delta P)^{1/2} - P^{1/2}}{\varepsilon}$$

Employing the chain rule, it is easy to see that

$$\delta S(P, \delta P)P^{1/2} + P^{1/2}\delta S(P, \delta P) = \delta P$$

so that

$$\delta S(P, \delta P) = \int_0^\infty \exp(-P^{1/2}t)\delta P \exp(-P^{1/2}t)\mathrm{d}t. \tag{45}$$

Taking (45) into account, we get the optimality condition

$$\int_0^\infty \left[\exp(-\hat{\Phi}_H^{1/2}t)\left(\hat{\Phi}_H^{1/2} - \Psi^{1/2}\right)\exp(-\hat{\Phi}_H^{1/2}t)\right]dt + \frac{1}{2}G^*\Lambda G = 0. \tag{46}$$

The integral in (46) is the unique solution of the Lyapunov equation

$$\hat{\Phi}^{1/2}X + X\hat{\Phi}^{1/2} = \hat{\Phi}^{1/2} - \Psi^{1/2}. \tag{47}$$

Equations (46)-(47) now yield

$$-\frac{1}{2}\hat{\Phi}^{1/2}(G^*\Lambda G) - \frac{1}{2}(G^*\Lambda G)\hat{\Phi}^{1/2} = \hat{\Phi}^{1/2} - \psi^{1/2},$$

which in turn gives

$$\hat{\Phi}^{1/2}\left(I + G^*\Lambda G\right) + \left(I + G^*\Lambda G\right)\hat{\Phi}^{1/2} = 2\Psi^{1/2}. \tag{48}$$

Since $I + G^*\Lambda G > 0$ almost everywhere on $\mathbb{T}$, we finally get

$$\hat{\Phi}^{1/2} = 2\int_0^\infty \exp\left[-(I + G^*\Lambda G)t\right]\Psi^{1/2}\exp\left[-(I + G^*\Lambda G)t\right]dt. \tag{49}$$

The maximization of the *dual functional* $\Lambda \to L^H(\hat{\Phi}, \Lambda)$, however, appears quite problematic.

We show in the next section that, differently from the Kullback-Leibler case, it is possible to define a sensible Hellinger distance for matricial functions that leads to a full unraveling of the complexity of the optimization problem. This will be accomplished by connecting this problem to a most classical topic at the hearth of systems and control theory, namely the *spectral factorization problem*.

## VI. HELLINGER DISTANCE AND SPECTRAL FACTORIZATION

Let $F$ be a measurable function defined on the unit circle $\mathbb{T}$ and taking values in $\mathbb{C}^{m \times p}$. Then $F$ belongs to the Hilbert space $L^2_{m \times p}$ if it satisfies

$$\int \operatorname{tr} F F^* d\theta < \infty.$$

For $F, G$ in $L^2_{m \times p}$, the scalar product is defined by

$$\langle F, G \rangle_2 = \int \operatorname{tr} F G^* d\theta,$$

so that $\|F\|_2^2 = \int \operatorname{tr} F F^* d\theta$. Let $\Phi \in \mathcal{S}^{m \times m}_+(\mathbb{T})$. Then a measurable $\mathbb{C}^{m \times p}$-valued function $W$ is called a spectral factor of $\Phi$ if it satisfies

$$W(e^{i\theta})W(e^{i\theta})^* = \Phi(e^{i\theta}), \quad \text{a.e. on } \mathbb{T}.$$

Notice that necessarily $p \geq m$ and $W(e^{i\theta})$ is a.e. full row rank. Moreover, $W$ is bounded on $\mathbb{T}$, and therefore it belongs to $L^2_{m \times p}$. When $p = m$, $W^{-1}$ is also bounded and, consequently, $W^{-1} \in L^2_{m \times m}$. Any $\Phi \in \mathcal{S}^{m \times m}_+(\mathbb{T})$ satisfies the Szegö condition

$$\int_{-\pi}^{\pi} \log \det \Phi(e^{i\theta}) \frac{d\theta}{2\pi} > -\infty,$$

and admits therefore spectral factors $W$ in $H^2_{m \times m}$, namely the Hardy space of functions in $L^2_{m \times m}$ that posses an analytic extension in $|z| > 1$, see e.g. [47], [35].

Let $W_1$ and $W_2$ be spectral factors of the same $\Phi \in \mathcal{S}^{m \times m}_+(\mathbb{T})$ with $W_1$ square. Then trivially $U := W_1^{-1} W_2$ is a $m \times p$ *all-pass* function, i.e

$$U(e^{i\theta})U(e^{i\theta})^* = I, \quad \forall e^{i\theta} \in \mathbb{T}.$$

For $\Phi, \Psi \in \mathcal{S}^{m \times m}_+(\mathbb{T})$, consider the following function

$$\tilde{d}_H(\Phi, \Psi) = \left[ \inf \left\{ \|W_\Psi - W_\Phi\|_2^2 : W_\Psi, W_\Phi \in L^2_{m \times m}, \ W_\Psi W_\Psi^* = \Psi, \ W_\Phi W_\Phi^* = \Phi \right\} \right]^{1/2}. \quad (50)$$

*Theorem 6.1: The following facts hold true:*

1) *For any square spectral factor $\bar{W}_\Psi$ of $\Psi$, we have:*

$$\tilde{d}_H(\Phi, \Psi) = \left[\inf\left\{\|\bar{W}_\Psi - W_\Phi\|_2^2 : \ W_\Phi \in L_{m\times m}^2, \ W_\Phi W_\Phi^* = \Phi\right\}\right]^{1/2}. \quad (51)$$

2) *The infimum in the above equation is a minimum: Indeed the unique spectral factor of $\Phi$ minimizing (51) is given by*

$$\hat{W}_\Phi := \Phi^{1/2}\left(\Phi^{1/2}\Psi\Phi^{1/2}\right)^{-1/2}\Phi^{1/2}\bar{W}_\Psi.$$

3) $\tilde{d}_H$ *is a bona fide* distance function.

4) $\tilde{d}_H$ *coincides with the Hellinger distance in the scalar case.*

   *Proof:*

1) First of all, observe that, once fixed the spectral factor $\bar{W}_\Psi$, any square spectral factor $W_\Psi$ of $\Psi$ can be written as $W_\Psi = \bar{W}_\Psi U$, where $U$ is a $m \times m$ all-pass. Hence,

$$\int \operatorname{tr}(W_\Psi - W_\Phi)(W_\Psi - W_\Phi)^* d\theta = \int \operatorname{tr}(\bar{W}_\Psi - W_\Phi U^*)(\bar{W}_\Psi - W_\Phi U^*)^* d\theta.$$

Observe, moreover, that $W_\Phi U^*$ is a square spectral factor of $\Phi$, so that (51) holds.

2) To show that the infimum in (51) is a minimum, notice that (51) may be rewritten in the form

$$\tilde{d}_H(\Phi, \Psi)^2 = \inf\left\{\int \operatorname{tr}(\bar{W}_\Psi - \Phi^{1/2}V)(\bar{W}_\Psi - \Phi^{1/2}V)^* d\theta : \ V \in L_{m\times m}^\infty, \ VV^* = I\right\}. \quad (52)$$

We shall solve this problem by unconstrained minimization of the Lagrangian

$$L = \int \operatorname{tr}\left[(\bar{W}_\Psi - \Phi^{1/2}V)(\bar{W}_\Psi - \Phi^{1/2}V)^* + \Delta(VV^* - I)\right],$$

where $\Delta = \Delta^* > 0$. The first variation of the Lagrangian (at $V$ in direction $\delta V \in L_{m\times m}^\infty$) is

$$\delta L(V; \delta V) = \int \operatorname{tr}\left[(\Delta V - \Phi^{1/2}W_\Psi)\delta V^* + \delta V(\Delta V - \Phi^{1/2}W_\Psi)^*\right]$$

The second variation of the Lagrangian is

$$\delta^2 L(V; \delta V) = 2\int \operatorname{tr}\left[\Phi^{1/2}\delta V\delta V^*\Phi^{1/2} + \Delta^{1/2}\delta V\delta V^*\Delta^{1/2}\right].$$

Hence, $L$ is strictly convex and therefore $V$ is a minimizer of the unconstrained minimization problem if and only if:

$$\delta L(V; \delta V) = 0, \qquad \forall \ \delta V. \quad (53)$$

Condition (53) is clearly equivalent to $\Delta V - \Phi^{1/2}W_\Psi = 0$ or to

$$V = \Delta^{-1}\Phi^{1/2}W_\Psi.$$

Thus, if there exists $\Delta = \Delta^* > 0$ such that

$$VV^* = \Delta^{-1}\Phi^{1/2}\Psi\Phi^{1/2}\Delta^{-1} = I,$$

then $V$ minimizes (52). Such a $\Delta$ is readily seen to be given by

$$\Delta = [\Phi^{1/2}\Psi\Phi^{1/2}]^{1/2}.$$

In conclusion, the infimum in (51) is a minimum and

$$\hat{W}_\Phi = \Phi^{1/2}\hat{V} = \Phi^{1/2}[\Phi^{1/2}\Psi\Phi^{1/2}]^{-1/2}\Phi^{1/2}\bar{W}_\Psi$$

is the unique minimizer.

3) The distance properties of $\tilde{d}_H$ are easy to check: (i) Symmetry is an immediate consequence of the definition of $\tilde{d}_H$. (ii) It is clear that $\tilde{d}_H(\Phi, \Phi) = 0$. Conversely, if $\tilde{d}_H(\Phi, \Psi) = 0$, then $\Phi$ and $\Psi$ share a common spectral factor and are, therefore, the same spectral density. (iii) The triangular inequality is inherited by the definition of $\tilde{d}_H$ as the infimum of the $L_2$ distance among spectral factors. Thus, given $\Phi$, $\Psi$ and $\Upsilon$ and chosen an arbitrary square spectral factor $W_\Upsilon$ of $\Upsilon$, we have

$$\tilde{d}_H(\Phi, \Psi) = \inf_{W_\Phi, W_\Psi} \|W_\Phi - W_\Psi\|_2 \leq \inf_{W_\Phi, W_\Psi} [\|W_\Phi - W_\Upsilon\|_2 + \|W_\Psi - W_\Upsilon\|_2] =$$
$$= \inf_{W_\Phi} \|W_\Phi - W_\Upsilon\|_2 + \inf_{W_\Psi} \|W_\Psi - W_\Upsilon\|_2 = \tilde{d}_H(\Phi, \Upsilon) + \tilde{d}_H(\Psi, \Upsilon)$$

where the last equality is a consequence of point 1).

4) By choosing $\bar{W}_\Psi = \Psi^{1/2}$, it is immediate to check that in the scalar case $(m = 1)$ $\hat{V} \equiv 1$ and hence $\tilde{d}_H$ coincides with the Hellinger distance $d_H$.                                                                                    $\square$

## VII. $\tilde{d}_H$-OPTIMAL MULTIVARIABLE SPECTRUM APPROXIMATION

Theorem 6.1 shows that $\tilde{d}_H$ is a natural extension to the multivariable case of the Hellinger distance. The corresponding multivariable version of Problem 3.4 is the following:

*Problem 7.1:* Given $\Psi \in \mathcal{S}_+^{m \times m}(\mathbb{T})$, find $\hat{\Phi} \in \mathcal{S}_+^{m \times m}(\mathbb{T})$ that solves

$$\text{minimize } \tilde{d}_H^2(\Phi, \Psi) \tag{54}$$

$$\text{subject to } \int G\Phi G^* = I. \tag{55}$$

It is in this form that the optimization problem is amenable to the variational analysis even in multivariable version. Let

$$\mathcal{L}^H := \{\Lambda \in \mathbb{C}^{n \times n} | \Lambda = \Lambda^*, I + G^*\Lambda G > 0 \ \forall e^{i\theta} \in \mathbb{T}\}, \tag{56}$$

and

$$\mathcal{L}_\Gamma^H := \mathcal{L}^H \cap \text{Range}(\Gamma) \tag{57}$$

where $\Gamma$ was defined in (28). The following is our main result.

*Theorem 7.2:* Assume condition (6) (or, equivalently, condition (8)) is satisfied. Then there exists a unique $\hat{\Lambda} \in \mathcal{L}_\Gamma^H$ such that

$$\int G(I + G^*\hat{\Lambda}G)^{-1}\Psi(I + G^*\hat{\Lambda}G)^{-1}G^* = I. \tag{58}$$

*The unique solution of the constrained approximation Problem 7.1 is then given by*

$$\hat{\Phi}_H := (I + G^*\hat{\Lambda}G)^{-1}\Psi(I + G^*\hat{\Lambda}G)^{-1}. \tag{59}$$

We break the proof of this theorem into two parts: First, by unconstrained minimization of the Lagrangian function, we obtain an expression for a spectral factor of the optimal $\Phi$ depending on the Lagrange multiplier matrix $\Lambda$ (Lemma7.3). Second, we establish existence of a unique $\Lambda \in \mathcal{L}_\Gamma^H$ satisfying (58) (Theorem7.6).

For $\Lambda \in \mathcal{L}^H$, $W_\Psi$ a spectral factor of $\Psi$, and $W, W^{-1} \in L_\infty^{m\times m}(\mathbb{T})$, form the Lagrangian function:

$$L(W, \Lambda) = \text{tr} \int (W - W_\Psi)(W - W_\Psi)^* + \text{tr} \, \Lambda \left( \int GWW^*G^* - I \right). \tag{60}$$

Consider the unconstrained minimization problem:

$$\min_W \left\{ L(W, \Lambda) | W, W^{-1} \in L_\infty^{m\times m}(\mathbb{T}) \right\} \tag{61}$$

*Lemma 7.3:* The unique solution to problem (61) is given by

$$\hat{W} = (I + G^*\Lambda G)^{-1}W_\Psi. \tag{62}$$

*Proof:* Let $\delta W \in L_\infty^{m\times m}(\mathbb{T})$. The first variation of the Lagrangian is:

$$\delta L(W, \Lambda; \delta W) = \text{tr} \int \left[ \delta W(W - W_\Psi)^* + (W - W_\Psi)\delta W^* + \Lambda(G\delta WW^*G^* + GW\delta W^*G^*) \right]$$

$$= \text{tr} \int (W - W_\Psi + G^*\Lambda GW)\delta W^* + \left( \text{tr} \int (W - W_\Psi + G^*\Lambda GW)\delta W^* \right)^*$$

By taking into account the cyclic property of the trace operator, the second variation of the Lagrangian is easily seen to be given by

$$\delta^2 L(W, \Lambda; \delta W) = 2\text{tr} \int \delta W^*(I + G^*\Lambda G)\delta W \tag{63}$$

which is clearly positive for any $\Lambda \in \mathcal{L}^H$ and $\delta W \neq 0$. Hence $L$ is strictly convex with respect to $W$. Moreover, the set $\mathcal{L}^H$ is open and convex. To find the minimum point of $L$, we impose $\delta L(W, \Lambda; \delta W) = 0$ in each direction $\delta W$. This yields (62). $\qquad \square$

We now consider the question of existence of a matrix $\hat{\Lambda} \in \mathcal{L}^H$ satisfying (58). To this end, we introduce the dual functional

$$
\begin{aligned}
L(\hat{W}, \Lambda) &= \operatorname{tr} \int \left( (I + G^*\Lambda G)^{-1} W_\Psi - W_\Psi \right) \left( (I + G^*\Lambda G)^{-1} W_\Psi - W_\Psi \right)^* \\
&\quad + \operatorname{tr} \left[ \Lambda \left( \int G (I + G^*\Lambda G)^{-1} W_\Psi W_\Psi^* (I + G^*\Lambda G)^{-1} G^* - I \right) \right] \\
&= \operatorname{tr} \int \Psi - (I + G^*\Lambda G)^{-1} \Psi - \operatorname{tr} \Lambda, \qquad \Lambda \in \mathcal{L}^H.
\end{aligned}
\tag{64}
$$

Instead of maximizing (64), we consider the equivalent problem of minimizing the functional:

$$
J_\Psi(\Lambda) := -L(\hat{W}, \Lambda) + \operatorname{tr} \int \Psi = \operatorname{tr} \int (I + G^*\Lambda G)^{-1} \Psi + \operatorname{tr} \Lambda, \qquad \Lambda \in \mathcal{L}^H.
\tag{65}
$$

*Lemma 7.4:* The functional (65) is convex and its restriction to $\mathcal{L}_\Gamma^H$ (defined in (57)) is strictly convex.

*Proof:* First of all, observe that $\mathcal{L}^H$ is an open, convex subset of the Hermitian matrices in $\mathbb{C}^{n \times n}$. For $\delta\Lambda \in \mathbb{C}^{n \times n}$ Hermitian, we compute the directional derivative:

$$
\begin{aligned}
\delta J_\Psi(\Lambda; \delta\Lambda) &= -\operatorname{tr} \int (I + G^*\Lambda G)^{-1} G^* \delta\Lambda G (I + G^*\Lambda G)^{-1} \Psi + \operatorname{tr} \delta\Lambda \\
&= \operatorname{tr} \left[ \left( I - \int G (I + G^*\Lambda G)^{-1} \Psi (I + G^*\Lambda G)^{-1} G^* \right) \delta\Lambda \right]
\end{aligned}
\tag{66}
$$

The second variation is then given by

$$
\delta^2 J_\Psi(\Lambda; \delta\Lambda) = 2\operatorname{tr} \int W_\Psi^* (I + G^*\Lambda G)^{-1} G^* \delta\Lambda G (I + G^*\Lambda G)^{-1} G^* \delta\Lambda G (I + G^*\Lambda G)^{-1} W_\Psi
\tag{67}
$$

which is clearly a non negative quantity. Hence $J_\Psi^H$ is convex on $\mathcal{L}^H$. In view of (31), we have that $\delta^2 J_\Psi(\Lambda; \delta\Lambda)$ is strictly positive for any non-zero $\delta\Lambda \in \operatorname{Range}(\Gamma)$, and consequently $J_\Psi$ is strictly convex on $\mathcal{L}_\Gamma^H$. $\qquad \square$

As an immediate consequence of the above lemma, we have the following corollary.

*Corollary 7.5:* The dual problem

$$
\text{Find } \Lambda \in \mathcal{L}_\Gamma^H \text{ minimizing } J_\Psi(\Lambda)
\tag{68}
$$

*admits at most one solution. Moreover, (58) is necessary and sufficient for $\hat{\Lambda}$ to solve the dual problem (68).*

We now tackle the existence issue for the dual problem. Although this is a finite-dimensional, convex optimization problem, the existence question is quite delicate since the set $\mathcal{L}_\Gamma^H$ is *open* and *unbounded*. The proof of the following theorem is partially inspired by the proof of the corresponding result for the scalar, Kullback-Leibler case in [Section 2][12].

*Theorem 7.6: If Problem 7.1 is feasible, i.e. (6) (or, equivalently, condition (8)) is satisfied, then the dual functional (65) has a unique minimum point in $\mathcal{L}_\Gamma^H$.*

*Proof:* In view of Corollary 7.5, we only need to show that $J_\Psi$ takes a minimum value on $\mathcal{L}_\Gamma^H$. First, we observe that $J_\Psi$ is continuous on its domain. Second, we show that $J_\Psi$ is bounded below on $\mathcal{L}_\Gamma^H$. Indeed, by feasibility, there exists a $\bar{\Phi} \in \mathcal{S}_+^{m \times m}(\mathbb{T})$ such that $\int G\bar{\Phi}G^* = I$. Hence, for all $M \in \mathbb{C}^{n \times n}$, $\int G\bar{\Phi}G^*M = M$, which implies

$$\operatorname{tr} M = \operatorname{tr} \int \bar{\Phi}^{1/2} G^* M G \bar{\Phi}^{1/2}. \tag{69}$$

Recalling that, for $\Lambda \in \mathcal{L}_\Gamma^H$, $I + G^*(e^{i\vartheta})\Lambda G(e^{i\vartheta})$ is positive definite for all $\vartheta \in [0, 2\pi)$, and using the monotonicity property of the trace, we get

$$\operatorname{tr} \Lambda = \operatorname{tr} \int \bar{\Phi}^{1/2} G^* \Lambda G \bar{\Phi}^{1/2} > -\operatorname{tr} \int \bar{\Phi}, \quad \forall \Lambda \in \mathcal{L}_\Gamma^H. \tag{70}$$

Define $\bar{f} := -\operatorname{tr} \int \bar{\Phi} < 0$. We get

$$J_\Psi(\Lambda) := \operatorname{tr} \int (I + G^*\Lambda G)^{-1}\Psi + \operatorname{tr} \Lambda = \operatorname{tr} \int \Psi^{1/2}(I + G^*\Lambda G)^{-1}\Psi^{1/2} + \operatorname{tr} \Lambda > \bar{f}, \qquad \forall \Lambda \in \mathcal{L}_\Gamma^H, \tag{71}$$

where we have used $\operatorname{tr} \int \Psi^{1/2}(I + G^*\Lambda G)^{-1}\Psi^{1/2} > 0.$ on $\mathcal{L}_\Gamma^H$.

Finally, we show that $J_\Psi$ is *inf-compact* i.e. the sub-level sets $J_\Psi^{-1}(-\infty, r]$ are compact. This implies existence of a minimum point. Indeed, observing that $J_\Psi(0) = \operatorname{tr} \int \Psi$, we can then restrict the search for a minimum point to the compact set $J_\Psi^{-1}(-\infty, \operatorname{tr} \int \Psi]$. Existence for the latter problem then follows from Weierstrass Theorem. To prove inf-compactness of $J_\Psi$, we proceed to show that:

  1)
$$\lim_{\Lambda \to \partial \mathcal{L}_\Gamma^H} J_\Psi(\Lambda) = +\infty,$$

  where $\partial \mathcal{L}_\Gamma^H$ denotes the boundary of $\mathcal{L}_\Gamma^H$;

  2)
$$\lim_{\|\Lambda\| \to \infty} J_\Psi(\Lambda) = +\infty.$$

To prove property 1), notice that $\partial \mathcal{L}_\Gamma^H$ is the set of $\Lambda$ in $\text{Range}(\Gamma)$ for which: (i) $I + G^* \Lambda G$ is positive semidefinite on $\mathbb{T}$ and (ii) $\exists \bar{\vartheta}$ s.t. $I + G^*(e^{i\bar{\vartheta}}) \Lambda G(e^{i\bar{\vartheta}})$ is singular. Thus, for $\Lambda \to \partial \mathcal{L}_\Gamma^H$, all the eigenvalues of $[I + G^* \Lambda G]^{-1}$ are positive on $\mathbb{T}$ and, at least one of them, has a pole tending to the unit circle ($[I + G^* \Lambda G]$ and $[I + G^* \Lambda G]^{-1}$ are *rational* matrix functions!). Since $\Psi$ is fixed and coercive, then also $\Psi^{1/2}[I + G^* \Lambda G]^{-1}\Psi^{1/2}$ has all eigenvalues positive on $\mathbb{T}$ and, at least one of them with a pole tending to the unit circle as $\Lambda \to \partial \mathcal{L}_\Gamma^H$. The conclusion follows by rewriting $J_\Psi(\Lambda)$ as $J_\Psi(\Lambda) = \text{tr} \int \Psi^{1/2}[I + G^* \Lambda G]^{-1}\Psi^{1/2} + \text{tr}\,\Lambda$, in view of (70).

Point 2) is more delicate. Let $\Lambda_k \in \mathcal{L}_\Gamma^H$ be a sequence such that $\lim_{k\to\infty} \|\Lambda_k\| = \infty$. Let $\Lambda_k^0 := \frac{\Lambda_k}{\|\Lambda_k\|}$. It is easy to see that if $\Lambda \in \mathcal{L}_\Gamma^H$, then $\alpha\Lambda \in \mathcal{L}_\Gamma^H$ for all $\alpha \in [0,1]$. Hence, for sufficiently large $k$, we have $\Lambda_k^0 \in \mathcal{L}_\Gamma^H$.

Let $\eta = \liminf \text{tr}\,\Lambda_k^0$. We want to show that $\eta$ is strictly positive. We first observe that $\eta \geq 0$. In fact, $\text{tr}\,\Lambda_k^0 = \frac{1}{\|\Lambda_k\|}\text{tr}\,\Lambda_k > \frac{1}{\|\Lambda_k\|}\bar{f} \to 0$, where we have used (70).

Consider a sub-sequence of $\Lambda_k^0$ such that the limit of its trace is $\eta$. Since this sub-sequence remains on the surface of the unit ball $\partial \mathcal{B} := \{\Lambda = \Lambda^* : \|\Lambda\| = 1\}$, which is compact, it has a sub-sub-sequence converging in $\partial \mathcal{B}$. Let $\Lambda_{k_i}^0$ be such a sub-sub-sequence, and let $\Lambda^\infty \in \partial \mathcal{B}$ be its limit. Clearly,

$$\lim_{i\to\infty} \Lambda_{k_i}^0 = \text{tr}\,\Lambda^\infty = \eta. \tag{72}$$

We now prove that $\Lambda^\infty \in \mathcal{L}_\Gamma^H$. To this aim, notice that $\Lambda^\infty$ is the limit of a sequence in the finite dimensional linear space $\text{Range}(\Gamma)$ and hence it belongs $\text{Range}(\Gamma)$. It remains to show that $(I + G^* \Lambda^\infty G)$ is positive definite on $\mathbb{T}$. Indeed, since (the unnormalized sub-sequence) $\Lambda_{k_i}$ belongs to $\mathcal{L}_\Gamma^H$, we have that $(I + G^* \Lambda_{k_i} G)$ is positive definite on $\mathbb{T}$ so that $(\frac{1}{\|\Lambda_{k_i}\|}I + G^* \Lambda_{k_i}^0 G)$ is also positive definite on $\mathbb{T}$. Hence, $G^* \Lambda^\infty G$ is positive semi-definite in $\mathbb{T}$ so that $(I + G^* \Lambda^\infty G)$ is strictly positive definite in $\mathbb{T}$. This proves that $\Lambda^\infty \in \mathcal{L}_\Gamma^H$. The latter, together with (69) yields

$$\text{tr}\,\Lambda^\infty = \text{tr} \int \bar{\Phi}^{1/2} G^* \Lambda^\infty G \bar{\Phi}^{1/2}. \tag{73}$$

As seen before, $G^* \Lambda^\infty G$ is positive semi-definite on $\mathbb{T}$. Moreover, $G^* \Lambda^\infty G$ is not identically zero since $\Lambda^\infty \in \text{Range}(\Gamma)$ (see (31)), and $\Lambda^\infty \neq 0$ (it is not the zero matrix) since $\Lambda^\infty \in \partial \mathcal{B}$. We conclude, in view of (72) and (73), that $\eta = \text{tr}\,\Lambda^\infty > 0$.

Finally, we have

$$J_\Psi(\Lambda_k) = \text{tr} \int \Psi^{1/2}(I + G^* \Lambda_k G)^{-1}\Psi^{1/2} + \text{tr}\,\Lambda_k \geq \|\Lambda_k\|\text{tr}\,\Lambda_k^0. \tag{74}$$

Since $\|\Lambda_k\| \to \infty$ and $\liminf \operatorname{tr} \Lambda_k^0 > 0$, we get

$$\lim_{k \to \infty} J_\Psi(\Lambda_k) = +\infty. \tag{75}$$

$\square$

Let $\hat{\Lambda} \in \mathcal{L}_\Gamma^H$ be the unique solution of the dual problem whose existence has just been proven in Theorem 7.6. We show below that it also provides via (59) the unique solution to the primal problem 7.1.

*Proof of Theorem 7.2:* Let $W_\Psi$ be any spectral factor of $\Psi$. Let $\hat{W} = (I + G^*\Lambda G)^{-1} W_\Psi$ as in (62). Let $W$, belonging to $L_\infty^{m \times m}(\mathbb{T})$ together with its inverse, satisfy the constraint

$$\int GWW^*G^* = I. \tag{76}$$

By Lemma (7.3), and by the strict convexity of the functional $L(\cdot, \hat{\Lambda})$, we get

$$\|\hat{W} - W_\Psi\|_2^2 = L(\hat{W}, \hat{\Lambda}) < L(W, \hat{\Lambda}) = \|W - W_\Psi\|_2^2.$$

Thus, $\hat{W}$ minimizes the $L_2$ distance to $W_\Psi$ among $W$ belonging to $L_\infty^{m \times m}(\mathbb{T})$ together with their inverse and satisfying constraint (76). Theorem 6.1 now shows that $\hat{\Phi}_H = \hat{W}\hat{W}^*$ (coinciding with $\hat{\Phi}_H$ in (59)), is the unique solution to the multivariate approximation Problem 7.1. $\square$

*Remark 7.7:* Consider the important *covariance extension* problem when, as it is often the case, the process $y$ is real-valued. Then $A$ and $B$ are real matrices and $\Psi$ is a real spectral density, i.e. $\Psi(z)$ is real (and symmetric) for all $z \in \mathbb{T}$. In this case, $\hat{\Lambda}$ is a real symmetric matrix.

## VIII. NUMERICAL SOLUTION OF THE DUAL PROBLEM

### A. A matricial Newton-type algorithm

We now show how to efficiently implement a *modified Newton algorithm* for the computation of $\hat{\Lambda}$ (convergence of the algorithm, however, will be discussed elsewhere). This task requires some care because we are working in a matricial space and vectorization does not appear to be convenient. The road-map is the following. We have to find the minimum of the functional (65) which is strictly convex on $\mathcal{L}_\Gamma^H$. This is then equivalent to finding a matrix $\hat{\Lambda} \in \mathcal{L}_\Gamma^H$ that annihilates the derivative of $J_\psi(\Lambda)$ i.e. such that (58) is satisfied. According to the abstract version of the Newton algorithm, $\hat{\Lambda}$ may be found as the limit of the sequence obtained by iterating the following steps:

1) Choose an initial estimate $\Lambda_0 \in \mathcal{L}_\Gamma^H$ of $\hat{\Lambda}$ (the simplest choice being $\Lambda_0 = 0$).

2) Let $\Lambda_i$ be the current estimate of $\hat{\Lambda}$. Compute the directional derivative $\delta J_\Psi(\Lambda_i; \delta\Lambda)$ at the point $\Lambda_i$ in direction $\delta\Lambda$ as in (66).

3) Compute the "Hessian" (second directional derivative) $\delta^2 J_\Psi(\Lambda_i; \delta\Lambda_1, \delta\Lambda_2)$ at the point $\Lambda_i$ in directions $\delta\Lambda_1$ and $\delta\Lambda_2$. This may be done in the same way in which we computed $\delta^2 J_\Psi(\Lambda_i; \delta\Lambda)$ in (67). Indeed, the latter may be viewed as the "diagonal" part of the Hessian i.e. $\delta^2 J_\Psi(\Lambda_i; \delta\Lambda) = \delta^2 J_\Psi(\Lambda_i; \delta\Lambda, \delta\Lambda)$. We get

$$\delta^2 J_\Psi(\Lambda_i; \delta\Lambda_1, \delta\Lambda_2) = \mathrm{tr} \int G Q_i^{-1} \left[ \left( G^* \delta\Lambda_2 G Q_i^{-1} \Psi \right) + \left( G^* \delta\Lambda_2 G Q_i^{-1} \Psi \right)^* \right] Q_i^{-1} G^* \delta\Lambda_1 \tag{77}$$

where $Q_i := (I + G^* \Lambda_i G)$.

4) Solve for $X \in \mathcal{L}_\Gamma^H$ the equation

$$\delta^2 J_\Psi(\Lambda_i; \delta\Lambda, X) = -\delta J_\Psi(\Lambda_i; \delta\Lambda) \qquad \forall \, \delta\Lambda \tag{78}$$

and set the $(i+1)$-th estimate of $\hat{\Lambda}$ to the value $\Lambda_{i+1} = \Lambda_i + X$

5) Let $\varepsilon$ be a suitably small number. If

$$\left\| \int G(I + G^* \Lambda_{i+1} G)^{-1} \Psi (I + G^* \Lambda_{i+1} G)^{-1} G^* - I \right\| > \varepsilon \tag{79}$$

then go to step 2). Otherwise, set $\hat{\Lambda} = \Lambda_{i+1}$.

There are some very delicate points to be addressed. First of all, we need to worry about existence of solutions for equation (78).

*Proposition 8.1:* Assume condition (6) (or, equivalently, condition (8)) is satisfied. There exists a unique $X \in \mathrm{Range}(\Gamma)$ solving Equation (78).

*Proof:* Equation (78) may be rewritten as

$$\int G Q_i^{-1} \left[ \left( G^* X G Q_i^{-1} \Psi \right) + \left( G^* X G Q_i^{-1} \Psi \right)^* \right] Q_i^{-1} G^* = \int G Q_i^{-1} \Psi Q_i^{-1} G^* - I \tag{80}$$

where we have eliminated $\delta\Lambda$. Notice that the map $\varphi$ associating to $X \in \mathrm{Range}(\Gamma)$ the matrix

$$\varphi(X) := \int G Q_i^{-1} \left[ \left( G^* X G Q_i^{-1} \Psi \right) + \left( G^* X G Q_i^{-1} \Psi \right)^* \right] Q_i^{-1} G^*$$

defines a linear transformation of $\mathrm{Range}(\Gamma)$ into itself. In fact, clearly

$$Q_i^{-1} \left[ \left( G^* X G Q_i^{-1} \Psi \right) + \left( G^* X G Q_i^{-1} \Psi \right)^* \right] Q_i^{-1} = \left[ Q_i^{-1} \left[ \left( G^* X G Q_i^{-1} \Psi \right) + \left( G^* X G Q_i^{-1} \Psi \right)^* \right] Q_i^{-1} \right]^* \tag{81}$$

so that by definition $\varphi(X) \in \mathrm{Range}(\Gamma)$. The linear map $\varphi$ has trivial kernel. In fact, if for some $X \in \mathrm{Range}(\Gamma)$, $\varphi(X) = 0$ then

$$0 = \mathrm{tr}\left[\varphi(X)X\right] = \delta^2 J_\Psi(\Lambda_i; X, X).$$

Taking into account the positive definiteness of $\delta^2 J_\Psi(\Lambda_i; X, X)$ on $\mathrm{Range}(\Gamma)$, this implies $X = 0$. As a consequence, the image of $\varphi$ is the whole linear space $\mathrm{Range}(\Gamma)$. It only remains to observe that $\left[\int G Q_i^{-1} \Psi Q_i^{-1} G^* - I\right] \in \mathrm{Range}(\Gamma)$. Indeed, $I \in \mathrm{Range}(\Gamma)$ and $\int G Q_i^{-1} \Psi Q_i^{-1} G^* \in \mathrm{Range}(\Gamma)$ by definition of $\mathrm{Range}(\Gamma)$. $\qquad\square$

Proposition 8.1 guarantees the existence of a solution of (78) (or, equivalently, of (80)) in $\mathrm{Range}(\Gamma)$ but not in $\mathcal{L}_\Gamma^H$ as requested. To overcome this problem, we resort to a variant of the Newton algorithm known as *Newton method with back-tracking*. In this variant, the following sub-steps are employed in place of step 4):

4.1) Solve for $X \in \mathrm{Range}(\Gamma)$ equation (78).

4.2) Let $k = 1$ and double $k$ until both conditions

$$\left(\Lambda_i + \frac{1}{k}X\right) \in \mathcal{L}_\Gamma^H, \tag{82}$$

$$\|V_i\| < \|V_{i,k}\|, \tag{83}$$

are satisfied, where

$$V_i := \int G(I + G^*\Lambda_i G)^{-1}\Psi(I + G^*\Lambda_i G)^{-1}G^* - I$$

and

$$V_{i,k} := \int G(I + G^*(\Lambda_i + \frac{1}{k}X)G)^{-1}\Psi(I + G^*(\Lambda_i + \frac{1}{k}X)G)^{-1}G^* - I.$$

4.3) Set the $(i+1)$-th estimate of $\hat{\Lambda}$ to the value $\Lambda_{i+1} = \Lambda_i + \frac{1}{k}X$

This procedure guarantees that each $\Lambda_i \in \mathcal{L}_\Gamma^H$ even when $X \notin \mathcal{L}_\Gamma^H$. Notice that, by convexity of the problem $X$ is a *descent direction* so that, for sufficiently large $k$, (83) is certainly satisfied. Moreover, since $\mathcal{L}_\Gamma^H$ is an open set, (82) is also satisfied for for sufficiently large $k$.

## B. Computation of the solution of equation (80)

The next point that needs to be addressed is the computation of $X$ (step 4.1)). In fact, although (80) is a linear equation, it is not obvious how to solve it in a numerically efficient way. To simplify notation, we drop the subscript "$i$" in $\Lambda_i$ and $Q_i := (I + G^*\Lambda_i G)$. Consider the following equation

$$\int GQ^{-1}\left[\left(G^*XGQ^{-1}\Psi\right) + \left(G^*XGQ^{-1}\Psi\right)^*\right]Q^{-1}G^* = \int GQ^{-1}\Psi Q^{-1}G^* - I \quad (84)$$

We propose the following procedure.

1) Choose a set $\{X_i\}$ of linearly independent matrices such that $\operatorname{span}\{X_i\} = \operatorname{Range}(\Gamma)$.

2) Compute the quantities

$$Y_i := \int GQ^{-1}\left[\left(G^*X_iGQ^{-1}\Psi\right) + \left(G^*X_iGQ^{-1}\Psi\right)^*\right]Q^{-1}G^*, \quad (85)$$

$$Y := \int GQ^{-1}\Psi Q^{-1}G^* - I. \quad (86)$$

3) Solve for the scalar unknowns $y_i$ equation

$$\sum_i y_i Y_i = Y.$$

4) Set

$$X = \sum_i y_i X_i.$$

Steps 3) and 4) do not present any difficulty. Concerning point 1), employing the characterization (29) of $\operatorname{Range}(\Gamma)$, we simply have to solve the $n \times m$ Lyapunov equations

$$\Sigma - A\Sigma A^* = BH_{h,k} + H_{h,k}^* B^* \quad (87)$$

where $H_{h,k}$ is the matrix in which the entry in position $(h,k)$ is 1 and all the other entries are zero. The $n \times m$ matrices obtained by solving (87) span $\operatorname{Range}(\Gamma)$, but are not linearly independent. It is easy, however, to employ the *singular value decomposition algorithm* (which is very stable and robust) and reduce to a basis $\{X_i\}$ of $\operatorname{Range}(\Gamma)$.

Concerning point 2), in the case when $\Psi$ is a *rational* matrix function we can compute the integrals in (85) and (86) very efficiently and precisely. We first describe in details the computation of $Y$. To compute $\int GQ^{-1}\Psi Q^{-1}G^*$ we observe that $\chi := GQ^{-1}\Psi Q^{-1}G^*$ is a spectral density. Let $W_\chi$ be an analytic spectral factor of $\chi$ (namely a function $W_\chi$ analytic outside the unit disk and such that $\chi = W_\chi W_\chi^*$). Then, $\int \chi$ is the steady-state covariance of the output of a filter with transfer function $W_\chi$ fed by normalized white noise. To compute a realization of $W_\chi$ we implement the following steps:

2a) Compute a *co-analytic* square spectral factor $W_\Psi^*$ of $\Psi$ (namely $\Psi = W_\Psi^* W_\Psi$, $W_\Psi$ being square and analytic outside the unit disk). This requires (see, e.g., [19]):

- To decompose $\Psi$ as $\Psi = Z + Z^*$ with $Z$ being analytic inside the unit disk: This may be done by partial fraction expansion.

- To solve an Algebraic Riccati equation of dimension equal to the McMillan degree of $Z$.

2b) Factorize $Q = (I + G^*\Lambda G)$ as $Q = \Delta^*\Delta$, with $\Delta$ being square and analytic together with its inverse outside the unit disk. This can be done by computing the stabilizing solution $X_s$ of the following algebraic Riccati equation:

$$X = A^*XA - A^*XB(I + B^*XB)^{-1}B^*XA + \Lambda \tag{88}$$

We have the following realization for $\Delta^{-1}$:

$$\Delta^{-1} = (I + B^*X_sB)^{-1/2} - (I + B^*X_sB)^{-1}B^*X_sA(zI - \Gamma_s)^{-1}B(I + B^*X_sB)^{-1/2} \tag{89}$$

with $\Gamma_s := A - B(I + B^*X_sB)^{-1}B^*X_sA$ having all its eigenvalues inside the unit circle.

2c) Compute a realization of $H^* := \Delta^{-*}W_\Psi^*$. Notice that $H^*$ is a co-analytic spectral factor of $\Delta^{-*}\Psi\Delta^{-1}$.

2d) From $H^*$, compute an analytic spectral factor $H_1$ of $\Delta^{-*}\Psi\Delta^{-1}$ using the procedure detailed in Lemma A.1 in the Appendix.

2e) Compute a minimal realization of $W_\chi := G\Delta H_1$:

$$W_\chi = C_\chi(zI - A_\chi)^{-1}B_\chi.$$

We get

$$\int GQ^{-1}\Psi Q^{-1}G^* = \int W_\chi W_\chi^* = C_\chi \Sigma_\chi C_\chi^*, \tag{90}$$

with $\Sigma_\chi$ being the solution of the Lyapunov equation

$$\Sigma_\chi - A_\chi \Sigma_\chi A_\chi^* = B_\chi B_\chi^*. \tag{91}$$

For the computation of the integrals $Y_i$ in (85), we employ the same technique. The main difference is that the integrand of (85) is not a spectral density. Nevertheless, we observe that, by factoring $Q$ as $Q = \Delta^*\Delta$ (exactly as we have done in point 2b) above) and by defining the functions $\Phi_1, \Phi_2 \in \mathcal{C}_H(\mathbb{T})$ as

$$\Phi_1 := \Delta^{-*}G^*XG\Delta, \qquad \Phi_2 := \Delta^{-*}\Psi\Delta, \tag{92}$$

we may rewrite such an integrand in the form

$$G\Delta^{-*}[\Phi_1\Phi_2 + \Phi_1\Phi_2]\Delta^{-*}G^* = G\Delta^{-*}[(\Phi_1 + \Phi_2)(\Phi_1 + \Phi_2)^* - \Phi_1\Phi_1^* - \Phi_2\Phi_2^*]\Delta^{-*}G^* \tag{93}$$

It is therefore clear that the integrand of (85) is a difference of spectral densities. Hence, the integral (85) may be computed by resorting to the same technique detailed above for the computation of $Y$.

## C. Simulation results

We have applied the procedure described in Section VIII-A to many different examples and it performed very well even in the case of large values of $n$ and $m$ (recall that $B \in \mathbb{C}^{n \times m}$). For example, for $n = 12$ and $m = 6$ the algorithm converged in less than 10 minutes in an Apple G4 1GHz computer. In the scalar case ($m = 1$) some examples are discussed in [20] for the case of the Hellinger distance and in [45] for the case of the KL pseudo-distance. In the following we discuss a simple multivariable example ($n = 3$, $m = 2$). Choose

$$A = \begin{bmatrix} 1/3 & 0 & 0 \\ 0 & 1/2 & 0 \\ 0 & 0 & 1/4 \end{bmatrix}, \qquad B = \begin{bmatrix} 1 & 0 \\ 1 & 1 \\ 0 & 1 \end{bmatrix}, \qquad \Sigma = \begin{bmatrix} 9/4 & 12/5 & 0 \\ 12/5 & 16/3 & 16/7 \\ 0 & 16/7 & 32/15 \end{bmatrix}.$$

$A$, $B$ and $\Sigma$ satisfy the feasibility condition (6). After re-normalizing so that $\Sigma = I$, we get:

$$A \simeq \begin{bmatrix} 0.2309 & -0.3657 & -0.2046 \\ 0.1368 & 0.7935 & 0.2434 \\ -0.0886 & -0.4086 & 0.0590 \end{bmatrix}, \qquad B \simeq \begin{bmatrix} 0.6191 & -0.0471 \\ 0.2697 & 0.2711 \\ -0.0458 & 0.6393 \end{bmatrix}$$

Finally, we have chosen the reference spectral density $\Psi$ to be identically equal to $I$ (the identity). In this case the Kullback-Leibler approximation has the interpretation of maximum entropy solution and may be obtained in closed form [26]:

$$\hat{\Phi}_{KL} = \left[ G^* B (B^* B)^{-1} B^* G \right]^{-1}. \tag{94}$$

For the computation of the Hellinger approximation, we have set in (79) $\varepsilon := \exp(-12)$. Our procedure converged in 12 steps of the Newotn algorithm with backtracking (less than 10 seconds) to the matrix

$$\hat{\Lambda}_H \simeq \begin{bmatrix} -0.4267 & -0.0621 & 0.0005 \\ -0.0621 & -0.1007 & -0.0269 \\ 0.0005 & -0.0269 & -0.4690 \end{bmatrix}.$$

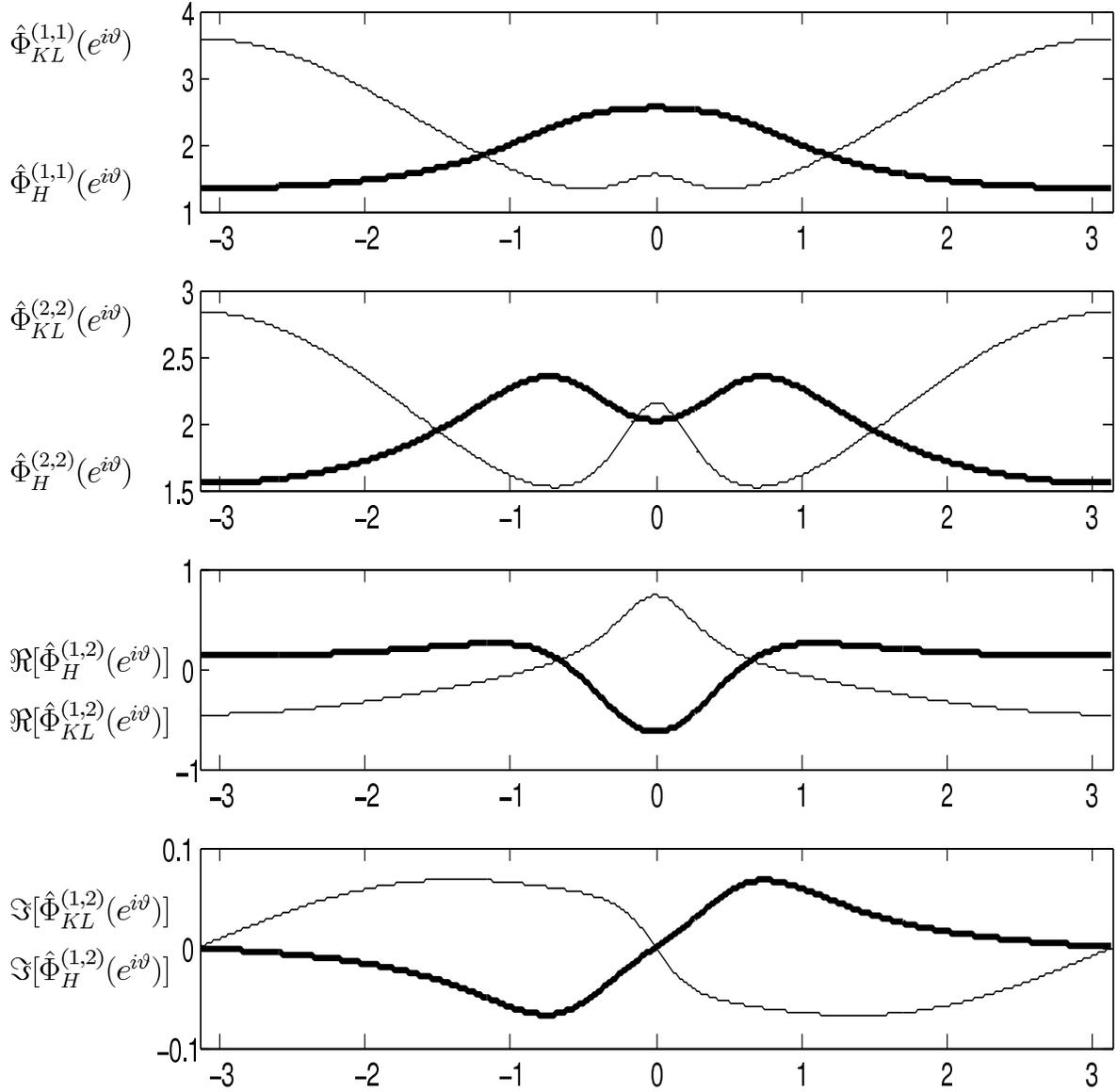

Fig. 1.   First picture: Graphics of $\hat{\Phi}_H^{(1,1)}(e^{i\vartheta})$ (bold line) and $\hat{\Phi}_{KL}^{(1,1)}(e^{i\vartheta})$ (thin line) as functions of $\vartheta$. Second picture: Graphics of $\hat{\Phi}_H^{(2,2)}(e^{i\vartheta})$ (bold line) and $\hat{\Phi}_{KL}^{(2,2)}(e^{i\vartheta})$ (thin line) as functions of $\vartheta$. Third picture: Graphics of $\Re[\hat{\Phi}_H^{(1,2)}(e^{i\vartheta})]$ (bold line) and $\Re[\hat{\Phi}_{KL}^{(1,2)}(e^{i\vartheta})]$ (thin line) as functions of $\vartheta$. Fourth picture: Graphics of $\Im[\hat{\Phi}_H^{(1,2)}(e^{i\vartheta})]$ (bold line) and $\Im[\hat{\Phi}_{KL}^{(1,2)}(e^{i\vartheta})]$ (thin line) as functions of $\vartheta$.

Let $\hat{\Phi}_H$ be the corresponding Hellinger approximation computed as in (59). Let $\hat{\Phi}_H^{(i,j)}$ be the entry in row $i$ and column $j$ of $\hat{\Phi}_H$ and similarly define $\hat{\Phi}_H^{(i,j)}$. In Figure 1, $\hat{\Phi}_H^{(1,1)}$, $\hat{\Phi}_H^{(2,2)}$ and the real and imaginary parts of $\hat{\Phi}_H^{(1,2)}$ are depicted together with the corresponding entries of $\hat{\Phi}_{KL}$.

It may be worthwhile to observe that, with respect to the $L_1$ distance (which induces a natural metric for spectral densities), Hellinger approximation does better than Kullback-Leibler approximation in this example. More precisely, we have

$$\mathrm{tr} \int \left[ \left( \hat{\Phi}_H - I \right)^2 \right]^{1/2} \simeq 1.6982$$

$$\mathrm{tr} \int \left[ \left( \hat{\Phi}_{KL} - I \right)^2 \right]^{1/2} \simeq 2.5079$$

## ACKNOWLEDGMENTS

The authors wish to thank prof. T. Georgiou for a very enlightening correspondence on approximation of spectral densities. They also wish to thank prof. A. Lindquist for illustrating to them various versions of the Byrnes-Lindquist existence theorem.

## APPENDIX

## A SPECTRAL FACTORIZATION RESULT

In this appendix we collect a side result on spectral factorization.

*Lemma A.1:* Let $A$ be a stability matrix and $H(z) = C(zI - A)^{-1}B + D$ be a minimal realization. Let $P$ be the solution of the Lyapunov equation

$$P - A^*PA = CC^*. \tag{95}$$

Let $\begin{bmatrix} K \\ J \end{bmatrix}$ be an ortho-normal basis of $\ker \begin{bmatrix} A^*P^{1/2} & C^* \end{bmatrix}$ i.e.

$$\begin{bmatrix} A^*P^{1/2} & C^* \end{bmatrix} \begin{bmatrix} K \\ J \end{bmatrix} = 0, \qquad \begin{bmatrix} K^* & J^* \end{bmatrix} \begin{bmatrix} K \\ J \end{bmatrix} = I. \tag{96}$$

Let $G := P^{-1/2}K$ and

$$H_1(z) := (D^*C + B^*PA)(zI - A)^{-1}G + B^*PG + D^*J. \tag{97}$$

Then, $H^*H = H_1H_1^*$.

*Proof:* Let $Q := C(zI - A)^{-1}G + J$. We have

$$Q^*Q = G^*(z^{-1}I - A^*)^{-1}C^*C(zI-A)^{-1}G + G^*(z^{-1}I - A^*)^{-1}C^*J + J^*C(zI-A)^{-1}G + J^*J \tag{98}$$

Now let $P > 0$ be the solution of the Lyapunov equation (95). Then,

$$C^*C = -(z^{-1}I - A^*)P(zI - A) + (z^{-1}I - A^*)Pz + z^{-1}P(zI - A) \tag{99}$$

Substituting (99) into (98) we obtain

$$
\begin{aligned}
Q^*Q &= -G^*PG + G^*Pz(zI-A)^{-1}G + G^*(z^{-1}I-A^*)^{-1}z^{-1}PG \\
&\quad + G^*(z^{-1}I-A^*)^{-1}C^*J + J^*C(zI-A)^{-1}G + J^*J
\end{aligned}
\tag{100}
$$

Moreover,

$$
z(zI-A)^{-1}I + A(zI-A)^{-1} \quad \text{and} \quad (z^{-1}I-A^*)^{-1}z^{-1} + I + (z^{-1}I-A^*)^{-1}A^*
\tag{101}
$$

so that

$$
Q^*Q = (J^*C + G^*PA)(zI-A)^{-1}G + \big((J^*C + G^*PA)(zI-A)^{-1}G\big)^* + G^*PG + J^*J
\tag{102}
$$

Taking (96) into account, it is easy to see that $Q^*Q = I$. Therefore, $H^*H = H^*QQ^*H$. Recalling (99) and (101), we eventually get

$$
\begin{aligned}
Q^*H &= (G^*(z^{-1}I-A^*)^{-1}C^* + J^*)(C(zI-A)^{-1}B + D) \\
&= -G^*(z^{-1}I-A^*)^{-1}(z^{-1}I-A^*)P(zI-A)(zI-A)^{-1}B \\
&\quad + G^*(z^{-1}I-A^*)^{-1}(z^{-1}I-A^*)Pz(zI-A)^{-1}B \\
&\quad + G^*(z^{-1}I-A^*)^{-1}z^{-1}P(zI-A)(zI-A)^{-1}B \\
&\quad + G^*(z^{-1}I-A^*)^{-1}C^*D + J^*C(zI-A)^{-1}B + J^*D \\
&= -G^*PB + G^*Pz(zI-A)^{-1}B + G^*(z^{-1}I-A^*)^{-1}z^{-1}PB \\
&\quad + G^*(z^{-1}I-A^*)^{-1}C^*D + J^*C(zI-A)^{-1}B + J^*D \\
&= -G^*PB + G^*P\left(I + A(zI-A)^{-1}\right)B + G^*\left(I + (z^{-1}I-A^*)^{-1}A^*\right)PB \\
&\quad + G^*(z^{-1}I-A^*)^{-1}C^*D + J^*C(zI-A)^{-1}B + J^*D \\
&= G^*PB + G^*PA(zI-A)^{-1}B + G^*(z^{-1}I-A^*)^{-1}A^*PB \\
&\quad + G^*(z^{-1}I-A^*)^{-1}C^*D + J^*C(zI-A)^{-1}B + J^*D \\
&= G^*(z^{-1}I-A^*)^{-1}(C^*D + A^*PB) + (G^*PA + J^*C)(zI-A)^{-1}B + G^*PB + J^*D \\
&= G^*(z^{-1}I-A^*)^{-1}(C^*D + A^*PB) + G^*PB + J^*D = H_1^*
\end{aligned}
\tag{103}
$$

$\square$